\documentclass[a4paper]{article} 
\usepackage[left=1.0cm,right=1.0cm,top=2.7cm]{geometry}
\usepackage{authblk}
\usepackage{CJKutf8}
\usepackage{yfonts}
\usepackage{lipsum}
\usepackage[utf8]{inputenc}
\title{Non-conformality of large deviations of moving average of the random walk in strongly mixing environment}

\author{Jiaming Chen\thanks{chen.jiaming@cims.nyu.edu}}
\affil{Courant Institute of Mathematical Sciences, New York University}
\date{\today}

\usepackage{graphicx}
\usepackage{amssymb}
\usepackage{amsmath}
\usepackage{amsthm}
\usepackage{tikz}
\usepackage{todonotes}
\usepackage{enumitem}
\usepackage{verbatim}

\usepackage[T2A]{fontenc}  
\newcommand{\CYRk}{\text{\char234}} 
\newcommand{\CYRu}{\text{\char244}}

\usepackage{centernot}
\usepackage{mathtools}
\usepackage{ stmaryrd }

\usepackage{hyperref}
\hypersetup{allcolors=black,
pdftitle={averaged law determines quenched law},
pdfpagemode=FullScreen}
\numberwithin{equation}{section}
\usepackage{bbm}

\usepackage{braket}
\usepackage{amsmath,amsthm,amssymb}
\usepackage{physics}
\usepackage{mathtools}
\usepackage{bm}
\usepackage{esvect}
\usepackage[english]{babel}
\usepackage{extarrows} 
\usepackage{mathrsfs}
\usepackage{esint}

\usepackage{setspace}
\usepackage{titlesec}
\titleformat{\subsection}[runin]
  {\normalfont\large\bfseries}{\thesubsection}{1em}{}

\usepackage{amsmath}

\numberwithin{equation}{section}
\usepackage{bbm}

\usepackage{braket}
\usepackage{amsmath,amsthm,amssymb}
\usepackage{physics}
\usepackage{mathtools}
\usepackage{bm}
\usepackage{esvect}
\usepackage[english]{babel}

\newtheorem{theorem}{Theorem}[section]

\theoremstyle{definition}

\theoremstyle{remark}

\usepackage[english]{babel}

\makeatletter
\renewenvironment{proof}[1][\proofname]{%
  \par\pushQED{\qed}\normalfont%
  \topsep6\p@\@plus6\p@\relax
  \trivlist\item[\hskip\labelsep\bfseries#1\@addpunct{.}]%
  \ignorespaces
}{%
  \popQED\endtrivlist\@endpefalse
}
\makeatother

\begin{document}
\maketitle

\begin{abstract}
    The quenched and annealed large deviations of the random walk in random environment is shown to conform on any compact set whenever the level of disorder is sufficiently low. In this work, we show that these two large deviations always disagree at some interior point of the natural domain of the random walk in strongly mixing environment, regardless of the level of disorder.
\end{abstract}



\section{Introduction}
Random walk in random environment (RWRE) is a fundamental model of statistical mechanics, describing motion of particles in highly disordered medium with random jump probabilities. This model was introduced from crystal growth \cite{Temkin} and turbulence of fluid via Lorentz gas description \cite{Sinai}. Various properties of this Markovian process has been investigated in general elliptic environments, including the transience and recurrence \cite{Alili} on the line as well as the existence of invariant measures \cite{Kozlov}. The renewal structure was first proposed in \cite{Kesten/Kozlov/Spitzer}. We also refer to \cite{Drewitz/Ramirez} for a distinguished review of recent progress on this model. Our main result in the following Theorem \ref{thm: 2.1} says that the quenched and annealed large deviations of the random walk in strongly mixing random environment always disagree at some interior point of the unit ball, regardless of the randomness of the environment, and thus extending Yilmaz \cite[Proposition 4]{Yilmaz} to dependent random field.

\section{Random walk in strongly mixing random environments}
   The random environment on $\mathbb{Z}^d$ is denoted as probability vectors $(\omega_x)_{x\in\mathbb{Z}^d}$ with $\sum_{e\in\mathbb{V}}\omega(x,e)=1$, $\forall~x\in\mathbb{Z}^d$. Here we use $\mathbb{V}\coloneqq\{e\in\mathbb{Z}^d:\,\abs{e}_1=1\}$. We use $\mathscr{F}_\Omega$ for the canonical product $\sigma$-algebra on $\mathbb{Z}^d$, and the exact distribution $(\mathbb{P},\mathscr{F}_\Omega)$, $\mathbb{E}\coloneqq E_{\mathbb{P}}$ will be specified later under different explicit mixing conditions. Throughout this paper, the environment $(\omega_x)_{x\in\mathbb{Z}^d}$ is assumed \textbf{uniformly elliptic}, i.e.~$\mathbb{P}(\omega(x,e)\geq\kappa)=1$ with any $x\in\mathbb{Z}^d$, $e\in\mathbb{V}$, for some absolute constant $0<\kappa<1/(2d)$.\par 
    Given any environmental law $(\mathbb{P},\mathscr{F}_\Omega)$, the disorder is defined as $\text{dis}(\mathbb{P})\coloneqq\inf\{\epsilon>0:\,\xi(x,e)\in[1-\epsilon,1+\epsilon],\,\mathbb{P}\text{-a.s.}\,\forall~x\in\mathbb{Z}^d,\,e\in\mathbb{V}\}$ with $\xi(x,e)\coloneqq\omega(x,e)/\mathbb{E}[\omega(x,e)]$ for each $x\in\mathbb{Z}^d$ and $e\in\mathbb{V}$. To descibe explicitly the mixing conditions on $(\mathbb{P},\mathscr{F}_\Omega)$, we first introduce some new notations. The distance $d_1(A,B)$ stands for the $\ell^1$-distance between sets $A,B\subseteq\mathbb{Z}^d$. For an absolute constant $r>0$, the law $(\mathbb{P},\mathscr{F}_\Omega)$ throughout this paper is assumed $r$-\textbf{Markovian}, i.e.
    \[
        \mathbb{P}\big((\omega_x)_{x\in V}\in\vdot\big|\mathscr{F}_{V^c}\big) = \mathbb{P}\big((\omega_x)_{x\in V}\in\vdot\big|\mathscr{F}_{\partial^r V}\big),\qquad\forall~\text{finite}~V\in\mathbb{Z}^d,\qquad\mathbb{P}\text{-a.s.,}
    \]
    where $\partial^rV\coloneqq\{z\in\mathbb{Z}^d\backslash V:\,\exists~y\in V~\text{s.t.}~\abs{z-y}_1\leq r\}$ and $\mathscr{F}_\Lambda\coloneqq\sigma(\omega_x,\,x\in\Lambda)$ for any $\Lambda\subseteq\mathbb{Z}^d$. Following the convention in \cite[Section 2]{Chen}, the environment $\mathbb{P}$ satisfies the strongly mixing condition $\textbf{(SMX)}_{C,g}$ if 
    \[
        \frac{d\mathbb{P}((\omega_x)_{x\in\Delta}\in\vdot|\eta)}{d\mathbb{P}((\omega_x)_{x\in\Delta}\in\vdot|\eta^\prime)} \leq \exp\bigg(C(\Delta,V^c)\sum_{x\in\Delta\cup\partial^r\Delta,y\in A\cup\partial^rA}e^{-g\abs{x-y}_1}\bigg),\qquad\forall~\text{finite}~\Delta\subseteq V\subseteq\mathbb{Z}^d\quad\text{with}~d_1(\Delta,V^c)\geq r
    \]
    and $A\subseteq V^c$, simultaneously for all pairs of configurations $\eta,\eta^\prime\in\Omega$ which agree on $V^c\backslash A$, $\mathbb{P}$-a.s., where $C(\Delta,V^c)>0$ if and only if $d_1(\Delta,V^c)< L_0$ for some absolute $L_0>r$.\par
    Having defined the first layer of our process, the random walk $(X_n)_{n\geq0}$ now travels on $\mathbb{Z}^d$ with jump probabilities given by $\omega$, i.e.~ for each $x\in\mathbb{Z}^d$ the law $P_{x,\omega}$ and $E_{x,\omega}\coloneqq E_{P_{x,\omega}}$ of this random walk starting from $x$ is prescribed by $P_{x,\omega}(X_0=x)=1$ and $P_{x,\omega}(X_{n+1}=y+e|X_n=y)=\omega(y,e)$ for all $y\in\mathbb{Z}^d$, $e\in\mathbb{V}$ and $n\in\mathbb{N}$. We call $P_{x,\omega}$ the quenched law of the RWRE. Averaging $P_{x,\omega}$ over $\omega$, we call the semi-direct product
    \[
        P_x( B)\coloneqq\int_{\Omega} P_{x,\omega}(B)\,d\mathbb{P(\omega)},\qquad\forall~B\in\mathscr{B}((\mathbb{Z}^d)^{\mathbb{N}})
    \]
    the annealed law of the RWRE. Here $\mathscr{B}((\mathbb{Z}^d)^{\mathbb{N}})$ stands for the Borel $\sigma$-algebra on $(\mathbb{Z}^d)^{\mathbb{N}}$.\par
    The nearest-neighbor trajectories imply that the velocities $|X_N/N|_1\leq1$ for all $N\geq1$. And thus we address $\mathbb{D}\coloneqq\{x\in\mathbb{R}^d:\,\abs{x}_1\leq1\}$ as the velocity surface and we denote $\partial\mathbb{D}_{d-2}\coloneqq\{x\in\partial\mathbb{D}:\,x_j=0\;\,\text{for some}~1\leq j\leq d\}$. In \cite[Theorem 2.1]{Chen}, when $\mathbb{P}$ satisfies $\textbf{(SMX)}_{C,g}$ the RWRE verifies the following quenched and annealed LDP on $\mathbb{Z}^d$ for any $d\geq1$:
    \begin{itemize}
        \item There exists the rate function $I_a:\mathbb{D}\to[0,\infty)$ which is convex and continuous on $\text{int}(\mathbb{D})$ so that for any Borel $G\subseteq\mathbb{R}^d$,
        \[
            -\inf_{x\in G^{\circ}} I_a(x)\leq \varliminf_{N\to\infty} \frac{1}{N}\log P_0\bigg(\frac{X_N}{N}\in A\bigg)\leq\varlimsup_{N\to\infty} \frac{1}{N}\log P_0\bigg(\frac{X_N}{N}\in A\bigg) \leq -\inf_{x\in\overline{G}} I_a(x).
        \]
        \item There exists the deterministic rate function $I_a:\mathbb{D}\to[0,\infty)$ which is convex and continuous on $\text{int}(\mathbb{D})$ such that $\mathbb{P}$-a.s. for any Borel $G\subseteq\mathbb{R}^d$,
        \[
            -\inf_{x\in G^{\circ}} I_q(x)\leq \varliminf_{N\to\infty} \frac{1}{N}\log P_{0,\omega}\bigg(\frac{X_N}{N}\in A\bigg)\leq\varlimsup_{N\to\infty} \frac{1}{N}\log P_{0,\omega}\bigg(\frac{X_N}{N}\in A\bigg) \leq -\inf_{x\in\overline{G}} I_q(x).
        \]
    \end{itemize}\
    Recent advances \cite{Chen} hava revealed how the level of disorder can force the equality between $I_a(\vdot)$ and $I_q(\vdot)$ on $\mathbb{D}$. When $\mathbb{P}$ satisfies $\textbf{(SMX)}_{C,g}$ the RWRE verifies the control on the threshold of disorder\footnote{Given any ergodic environment $\mathbb{P}$, the disorder $\text{dis}(\mathbb{P})\coloneqq\inf\{\epsilon>0:\,\frac{\omega(x,e)}{\mathbb{E}[\omega(x,e)]}\in[1-\epsilon,1+\epsilon],\,\forall~e\in\mathbb{V}\;\;\text{and}\;\;x\in\mathbb{Z}^d\}$ is defined as the level of randomness of the environment.} below which the two rate functions conform: 
    \begin{itemize}
        \item For any $d\geq4$ and compact $\mathcal{K}\subseteq\partial\mathbb{D}\backslash\partial\mathbb{D}_{d-2}$, there exists $\epsilon_1=\epsilon_1(\mathcal{K})>0$ such that $I_a(x)=I_q(x)$ for any $x\in\mathcal{K}$, whenever $\text{dis}(\mathbb{P})<\epsilon_1$.
        \item For any $d\geq4$ and compact $\mathcal{K}\subseteq\text{int}(\mathbb{D})\backslash\{0\}$, there exists $\epsilon_2=\epsilon_2(\mathcal{K})>0$ such that $I_a(x)=I_q(x)$ for any $x\in\mathcal{K}$, whenever $\text{dis}(\mathbb{P})<\epsilon_2$.
    \end{itemize}
    Although controlling the level of disorder allows us to enforce the equality $I_a=I_q$ uniformly on compact subsets of $\mathbb{D}$, it is striking that they can never collide simultaneously on the entire velocity surface whenever the environment $\mathbb{P}$ exhibits genuine randomness. This fundamental non-conformality persists even under strongly mixing condition $\textbf{(SMX)}_{C,g}$ and extends a classical result of Yilmaz \cite[Proposition 4]{Yilmaz}, originally established for i.i.d.\ environments
    \begin{theorem}\label{thm: 2.1}
        \normalfont
        For any $d\geq1$ and $\kappa>0$, whenever $\mathbb{P}$ verifies $\textbf{(SMX)}_{C,g}$ and whose support $\mathbb{P}$ is not a singleton, i.e.~$\text{dis}(\mathbb{P})>0$, then $I_a(x)<I_q(x)$ at some interior point $x\in\text{int}(\mathbb{D})$ of the velocity surface.
    \end{theorem}
    The non-conformality of two rate functions $I_a$ and $I_q$ occurs exactly because they cannot be equal at one point of the edge-vertex set $\partial\mathbb{D}_{d-2}$. However, this property remains open in the more general $\textbf{(SM)}_{C,g}$ or $\textbf{(SMG)}_{C,g}$ scenarios.

\section{Construction of auxiliary $Q^z$-walk and separation of correlated field}\label{sec: 3}
    Take an arbitrary $z\in\text{int}(\mathbb{D})\backslash\{0\}$, choose $C^{\mathbb{P}}_{z}>0$ which ensures $ f(C^{\mathbb{P}}_{z})\coloneqq \frac{1}{2}\sum_{e\in\mathbb{V}} \sqrt{\abs{\langle  z,e\rangle}^2+ 4 C^{\mathbb{P}}_{z} \mathbb{E}[\omega(0,e)]\mathbb{E}[\omega(0,-e)]}=1$. Such $0<C^{\mathbb{P}}_{z}<\infty$ exists because $f(0)=\abs{z}<1$ and $f(\infty)=\infty$. Define the probability vector $(u_\lambda(e))_{e\in\mathbb{V}}$ by 
    \[
        2u_z(e)\coloneqq \langle z,e\rangle + \sqrt{ \abs{\langle  z,e\rangle}^2+ 4 C^{\mathbb{P}}_{z} \mathbb{E}[\omega(0,e)]\mathbb{E}[\omega(0,-e)]},\qquad\forall~e\in\mathbb{V}.
    \]
    Define the $Q^z$-walk $(Z_n)_{n\geq0}$ on $\mathbb{Z}^d$ starting from any $x\in\mathbb{Z}^d$ with law $Q^z_x$ prescribed by $ Q^z_x(Z_{n+1}=Z_n+e) = u_z(e)$ for any $e\in\mathbb{V}$ and $n\geq0$. With this $(u_z(e))_{e\in\mathbb{V}}$, there exists $c_z=c_z(\kappa)>0$ such that $u_z(e)\geq c_z$ for all $e\in\mathbb{V}$. Also $E^{Q^z}_x[Z_{n+1}-Z_n]= z$ for all $n\in\mathbb{N}$. Proceeding in analogous steps of \cite[Lemma 3.2]{Bazaes/Mukherjee/Ramirez/Sagliett0}, for each $n\geq1$ and $\theta\in\mathbb{R}^d$ we have
    \begin{equation}\label{eqn: Q walk form1, annealed}    
        E^{Q^{z}}_0\bigg[ e^{\langle\theta,Z_n\rangle}\mathbb{E} \prod_{j=1}^n\xi(Z_{j-1},\Delta_j(Z)) \bigg] = \mathcal{D}_{z,\mathbb{P}}^n E_{0}\bigg[ e^{\langle\theta+\theta^{\mathbb{P}}_{z},X_n\rangle} \bigg],\quad E^{Q^{z}}_0\bigg[ e^{\langle\theta,Z_n\rangle}\prod_{j=1}^n\xi(Z_{j-1},\Delta_j(Z)) \bigg] = \mathcal{D}_{z,\mathbb{P}}^n E_{0,\omega}\bigg[ e^{\langle\theta+\theta^{\mathbb{P}}_{z},X_n\rangle} \bigg]
    \end{equation}
    where $\mathcal{D}_{z,\mathbb{P}}\coloneqq (C^{\mathbb{P}}_{z})^{1/2}$ and the vector $\theta^{\mathbb{P}}_{z}\in\mathbb{R}^d$ is defined via $ \langle \theta^{\mathbb{P}}_{z},e_j\rangle\coloneqq \text{log}(\mathcal{D}_{z,\mathbb{P}}^{-1}u_z(e_j)/\mathbb{E}[\omega(0,e_j)])$, $\forall~j=1,\ldots,d$. And $\Delta_j(Z)$ stands for $Z_j-Z_{j-1}$ with $j=1,\ldots,n$. Following \cite[Theorem 4.3.1]{Dembo/Zeitouni} for the annealed scenario and \cite[Theorem 2.6]{Rassoul-Agha/Seppalainen} for the quenched scenario, from the above (\ref{eqn: Q walk form1, annealed}) it yields the limiting identities 
    \[
        \overline{\Lambda}^a_z(\theta)\coloneqq\lim_{N\to\infty}\frac{1}{N}\log E^{Q^{z}}_0\bigg[ e^{\langle\theta,Z_N\rangle}\mathbb{E} \prod_{j=1}^N\xi(Z_{j-1},\Delta_j(Z)) \bigg],\;\; \overline{\Lambda}^q_{z}(\theta)\coloneqq\lim_{N\to\infty}\frac{1}{N}\log E^{Q^{z}}_0\bigg[ e^{\langle\theta,Z_N\rangle} \prod_{j=1}^N\xi(Z_{j-1},\Delta_j(Z)) \bigg],\;\;\forall~\theta\in\mathbb{R}^d
    \]
    and satisfy $\overline{\Lambda}^a_{z}(\theta) = \log(C^{\mathbb{P}}_{z,})^{1/2} + \Lambda_a(\theta+\theta^{\mathbb{P}}_{z})$, $ \overline{\Lambda}^q_{z}(\theta) = \log(C^{\mathbb{P}}_{z})^{1/2} + \Lambda_q(\theta+\theta^{\mathbb{P}}_{z})$. We also use $\overline{\Lambda}^*_z(\theta)\coloneqq\log \sum_{z\in\mathbb{V}}e^{\langle\theta,e\rangle}u_z(e)$ to stand for the limiting free energy in the zero-disorder scenario, $\forall~\theta\in\mathbb{N}$. Fix $\ell\in\mathbb{V}$ such that $\langle z,\ell\rangle>0$. Let $\mathcal{W}\coloneqq\mathbb{V}\cup\{0\}$. Consider the probability measure $\overline{Q}^z_0$ given by $\overline{Q}^z_0\coloneqq U\otimes\overline{Q}^z_{\epsilon,0}$ on $  (\mathcal{W})^{\mathbb{N}}\times(\mathbb{Z}^d)^{\mathbb{N}}$, where $U$ is product measure on $(\mathcal{W})^{\mathbb{N}}$ s.t. $U(\epsilon_1=e)=\CYRk$, $\forall~e\in\mathbb{V}$, $U(\epsilon_1=0)=1-2d\CYRk$, $\forall~\epsilon=(\epsilon_1,\epsilon_2,\ldots)\in(\mathcal{W})^{\mathbb{N}}$ for some small $\CYRk=\CYRk(\kappa)>0$, and the Markov chain $\overline{Q}^z_{\epsilon,0}$ on $\mathbb{Z}^d$ is defined via 
    \[
        \overline{Q}^z_{\epsilon,0}(Z_{n+1}=x+e|Z_n=x) = \mathbbm{1}_{\{ \epsilon_{n+1}=e \} } +\frac{\mathbbm{1}_{\{ \epsilon_{n+1}=0 \} }}{1-2d\CYRk}\big(u_z(e)-\CYRk\big),\qquad\forall~e\in\mathbb{V},\quad x\in\mathbb{Z}^d.
    \]
    So here we can actually take $2\CYRk(\kappa)\coloneqq d^{-1}\wedge\inf\{u_z(e):\,e\in\mathbb{V}\}$. A crucial remark is that the laws $\overline{Q}^z_0 = U\otimes\overline{Q}^z_{\epsilon,0}$ and $Q^z_0$ coincide. We also write the sequence $\Bar{\epsilon}\coloneqq\ell$ and let $\Bar{\epsilon}^{(L)}\coloneqq(\Bar{\epsilon},\ldots,\Bar{\epsilon})$ stand for the $L$-repetition of $\Bar{\epsilon}$, $\forall~L\geq L_0>r$ while $\mathbb{P}$ assumed $r$-Markovian. We introduce the $Q$-stopping sequence $(\tau^{(L)}_n)_{n\geq0}$ by $\tau^{(L)}_0=0$ and
    \[
        \tau^{(L)}_{n}\coloneqq\inf\{j\geq\tau^{(L)}_{n-1}+L:\,(\epsilon_{j-L},\ldots,\epsilon_{j-1}) = +s,\ldots,s\},\qquad s=\ell, \qquad\forall~n\geq1.
    \]
    Remark that between each $\tau^{(L)}_{j-1}$ and $\tau^{(L)}_{j}$ the $\epsilon$-sequence creates a spacing where $(Z_n)_{n\geq0}$ under the law $\overline{Q}^z_0$ travels with full probability at $
    \ell$-direction of length $L$. Accordingly we define the auxiliary transition
    \[
        \psi_k(\CYRu)\coloneqq \prod_{j=\tau^{(L)}_{k-1}+1}^{\tau^{(L)}_{k}} \CYRu_j\;\;\text{with}\;\, \CYRu_j\coloneqq \mathbbm{1}_{ \{\epsilon_j=\Delta_j(Z)\} } + \mathbbm{1}_{ \{\epsilon_j=0\} }\bigg(\frac{\omega(Z_{j-1},\Delta_j(Z))}{\mathbb{E}[\omega(Z_{j-1},\Delta_j(Z))]} + \frac{\CYRk}{u_z(\Delta_j(Z))-\CYRk}\bigg( \frac{\omega(Z_{j-1},\Delta_j(Z))}{\mathbb{E}[\omega(Z_{j-1},\Delta_j(Z))]} - 1 \bigg) \bigg)
    \]
    for each $k,j\in\mathbb{N}$. Here we write $\tau^{(L)}_{0}=0$ for convenience. This auxiliary transition renders us with $E^{\overline{Q}^z}_0[e^{\langle\theta,Z_n\rangle}\prod_{j=1}^n\CYRu_j] = E^{Q^z}_0[e^{\langle\theta,Z_n\rangle}\prod_{j=1}^n\xi(Z_{j-1},\Delta_j(Z))]$, $\forall~\theta\in\mathbb{R}^d$ and $n\geq1$.

\section{Non-uniform conformality of rate functions in the domain interior}
    Denote by $\mathcal{S}_N$ the set of all nearest-neighbor paths of length $N$ on $\mathbb{Z}^d$. Then $\forall~z\in\text{int}(\mathbb{D})\backslash\{0\}$ and $\ell\in\mathbb{V}$,
    \[
        E^{\overline{Q}^z}_0\bigg[ e^{-\langle\theta^{\mathbb{P}}_z,Z_N\rangle}\prod_{j=1}^N\CYRu_j I_{ \{ Z_j=\ell j \} } \bigg] = \mathcal{D}^N_{z,\mathbb{P}}P_{0,\omega}(X_N=\ell N),\quad E^{\overline{Q}^z}_0\bigg[ e^{-\langle\theta^{\mathbb{P}}_z,Z_N\rangle}\mathbb{E}\prod_{j=1}^N\CYRu_j I_{ \{ Z_j=\ell j \} } \bigg] = \mathcal{D}^N_{z,\mathbb{P}}P_{0}(X_N=\ell N),
    \]
    for any $N\geq1$. Define the stopping sequence $(\tau^{(L)}_n)_{n\geq1}$ as in Section \ref{sec: 3} with respect to the $U$-probability. Then,
   \begin{equation}\begin{aligned}\label{eqn: express Ia}  
        -I_a(\ell) &\geq \varlimsup_{N\to\infty} \frac{1}{N}\log P_0(X_N=\ell N) \geq -\text{log}( e^{\langle \theta_{z}^{\mathbb{P}},\ell\rangle} \mathcal{D}_{z,\mathbb{P}} ) + \varlimsup_{N\to\infty} \frac{1}{N}\log E^{\overline{Q}^z}_0\bigg[\mathbb{E}\prod_{j=1}^N \CYRu_j I_{ \{ Z_j=\ell j \} } \bigg]\\
        &\geq -\text{log}( e^{\langle \theta_{z}^{\mathbb{P}},\ell\rangle} \mathcal{D}_{z,\mathbb{P}} ) + \varlimsup_{k\to\infty} \frac{1}{E^{\overline{Q}^z}_0[\tau^{(L)}_k]}\log E^{\overline{Q}^z}_0 \bigg[\mathbb{E}\prod_{j=1}^{\tau^{(L)}_k} \CYRu_j I_{ \{ Z_j=\ell j \} } \bigg].
   \end{aligned}\end{equation}
   where the last line follows from Lemma \cite[Lemma 4.4]{Chen}. In the quenched scenario, $\forall~\epsilon>0$ we can write
    \begin{equation*}\begin{aligned}
        -I_q(\ell) &\leq \varliminf_{N\to\infty} \frac{1}{N}\log P_{0,\omega}(\langle X_N,\ell\rangle>(1-\epsilon/E^{\overline{Q}^z}_0[\tau^{(L)}_1])N)\leq -\text{log}( e^{\langle \theta_{z}^{\mathbb{P}},\ell\rangle} \mathcal{D}_{z,\mathbb{P}} )\\
        &\quad+ \varliminf_{N\to\infty} \frac{1}{N}\log E^{\overline{Q}^z}_0 \bigg[\prod_{j=1}^N \CYRu_j,\,\langle Z_N,\ell\rangle>(1-\epsilon/E^{\overline{Q}^z}_0[\tau^{(L)}_1])N \bigg] +\epsilon\log e^{\langle \theta_{z}^{\mathbb{P}},\ell\rangle}.
    \end{aligned}\end{equation*}
    And therefore we derive the estimate
    \begin{equation*}\begin{aligned}
        -I_q(\ell) &\leq -\text{log}( e^{\langle \theta_{z}^{\mathbb{P}},\ell\rangle} \mathcal{D}_{z,\mathbb{P}} ) + \varliminf_{k\to\infty} \frac{1}{E^{\overline{Q}^z}_0[\tau^{(L)}_k]}\log E^{\overline{Q}^z}_0\bigg[\prod_{j=1}^{\tau^{(L)}_k} \CYRu_j,\,\langle Z(\tau^{(L)}_k),\ell\rangle> (1-\epsilon/E^{\overline{Q}^z}_0[\tau^{(L)}_1])\tau^{(L)}_k  \bigg]  + \epsilon\log e^{\langle \theta_{z}^{\mathbb{P}},\ell\rangle}\\
        &\leq -\text{log}( e^{\langle \theta_{z}^{\mathbb{P}},\ell\rangle} \mathcal{D}_{z,\mathbb{P}} ) + (1-\epsilon/E^{\overline{Q}^z}_0[\tau^{(L)}_1]) \varliminf_{k\to\infty} \frac{1}{E^{\overline{Q}^z}_0[\tau^{(L)}_k]}\log E^{\overline{Q}^z}_0\bigg[\prod_{j=1}^{\tau^{(L)}_k} \CYRu_j I_{ \{ Z_j=\ell j \} } \bigg] + \epsilon\log e^{\langle \theta_{z}^{\mathbb{P}},\ell\rangle} + O(\epsilon),
    \end{aligned}\end{equation*}
    since for each $N\geq1$, the number of paths constituting the event $\{\langle Z_N,\ell\rangle>(1-\epsilon/E^{\overline{Q}^z}_0[\tau^{(L)}_1])N\}$ is $e^{ O(\epsilon)N}$, where there are at least $(1-\epsilon/E^{\overline{Q}^z}_0[\tau^{(L)}_1])N$ such jumps advancing in the $\ell$-direction. Let $\epsilon\to0$,
    \begin{equation}\label{eqn: express Iq}
        -I_q(\ell) \leq -\text{log}( e^{\langle \theta_{z}^{\mathbb{P}},\ell\rangle} \mathcal{D}_{z,\mathbb{P}} ) + \frac{1}{E^{\overline{Q}^z}_0[\tau^{(L)}_1]}\mathbb{E}\log E^{\overline{Q}^z}_0\bigg[\prod_{j=1}^{\tau^{(L)}_1} \CYRu_j I_{ \{ Z_j=\ell j \} } \bigg].
    \end{equation}
    \begin{proof}[\textbf{Proof of Theorem \ref{thm: 2.1}}.]
        Invoking a separation lemma for the stopping sequence $(\tau^{(L)}_k)_{k\geq1}$ via \cite[Lemma 4.1]{Chen}, we have
        \begin{equation}\label{eqn: express Ia, sharp}
            -I_a(\ell) \geq -\text{log}( e^{\langle \theta_{z}^{\mathbb{P}},\ell\rangle} \mathcal{D}_{z,\mathbb{P}} ) +   \frac{1}{E^{\overline{Q}^z}_0[\tau^{(L)}_1]}\log E^{\overline{Q}^z}_0\bigg[\mathbb{E}\prod_{j=1}^{\tau^{(L)}_1} \CYRu_jI_{ \{ Z_j=\ell j\} } \bigg].
        \end{equation}
        Writing $W^z_{\mathbb{P}}\coloneqq\text{log}( e^{\langle \theta_{z}^{\mathbb{P}},\ell\rangle} \mathcal{D}_{z,\mathbb{P}} )$ and combining (\ref{eqn: express Iq}) with (\ref{eqn: express Ia, sharp}), we have       
    \[
        \frac{W^z_{\mathbb{P}}-I_q(\ell)}{W^z_{\mathbb{P}}-I_a(\ell)} \leq \frac{\mathbb{E}\log E^{\overline{Q}^z}_0 [\prod_{j=1}^{T^{(L)}_1}  \CYRu_jI_{ \{ Z_j=\ell j\} } ]}{\log E^{\overline{Q}^z}_0[\mathbb{E}\prod_{j=1}^{T^{(L)}_1}  \CYRu_jI_{ \{ Z_j=\ell j\} } ] }.
    \]
    Notice also that the support of $\mathbb{P}$ is not a singleton for each $x\in\mathbb{Z}^d$ and $\ell\in\mathbb{V}$, i.e.~we have $\mathbb{P}(\omega(x,\ell)=\mathbb{E}[\omega(x,\ell)])<1$. Henceforth, by Jensen's inequality, 
    \begin{equation}\begin{aligned}\label{eqn: what we guess}
        \frac{1}{E^{\overline{Q}^z}_0[\tau^{(L)}_1]} \mathbb{E} \log E^{\overline{Q}^z}_0\bigg[\prod_{j=1}^{\tau^{(L)}_1} \CYRu_jI_{ \{ Z_j=\ell j\} } \bigg] < \frac{1}{E^{\overline{Q}^z}_0[\tau^{(L)}_1]}\log E^{\overline{Q}^z}_0\bigg[\mathbb{E}\prod_{j=1}^{\tau^{(L)}_1} \CYRu_jI_{ \{ Z_j=\ell j\} } \bigg],
    \end{aligned}\end{equation}
    yielding that $I_q(\ell)>I_a(\ell)$ with $\ell\in\mathbb{V}$. Since both $I_a$ and $I_q$ are continuous on the velocity surface $\mathbb{D}$, we can find some $x\in\text{int}(\mathbb{D})$ close to $\ell\in\mathbb{V}$ such that $I_q(x)>I_a(x)$, verifying the assertion.
    \end{proof}

\appendix

\bibliographystyle{plain}
\begin{spacing}{1}

\end{spacing}

\end{document}